\documentclass[a4 paper,twoside]{amsart}

\usepackage{graphics}
\usepackage{latexsym}
\usepackage{amsfonts}
\usepackage{amssymb}

\def\RR{\mathbb{R}}
\def\MM{\mathbb{M}}

\def\dto{{\scriptstyle\buildrel{\scriptstyle\longrightarrow}
\over{{}_{\scriptstyle\longrightarrow}}}}
\def\eps{\varepsilon}
\def\wto{\rightharpoonup}
\def\union{\mathop\cup}
\def\dom{{\rm dom}}
\def\Z{\mathcal{Z}}
\def\det{{\rm det}}
\def\I{E}
\def\J{\mathcal{\I}}
\def\Aff{{\rm Aff}}
\def\Cup{\mathop\cup}
\def\Cap{\mathop\cap}
\def\E{\J}

\def\Z{\mathcal{Z}}
\def\JJ{\mathcal{I}}
\def\Omega{\Sigma}

\newtheorem{theorem}{Theorem}[section]
\newtheorem{lemma}[theorem]{Lemma}
\newtheorem{proposition}[theorem]{Proposition}
\newtheorem{corollary}[theorem]{Corollary}

\theoremstyle{definition}
\newtheorem{definition}[theorem]{Definition}

\theoremstyle{remark}
\newtheorem{remark}[theorem]{Remark}

\title[The nonlinear membrane energy]{The nonlinear membrane energy: variational derivation under the constraint \boldmath$``\det\nabla u\not=0"$\unboldmath}

\author{Omar Anza Hafsa}
\address{Institute of Mathematics, University of Z\"urich, Winterthurerstrasse 190, CH-8057 Z\"urich, Switzerland.}
\email{anza@math.unizh.ch}

\author{Jean-Philippe Mandallena}
\address{``Equipe AVA (Analyse Variationnelle et Applications)" Centre Universitaire de Formation et de Recherche de N\^\i mes, Site des Carmes, Place Gabriel P\'eri - Cedex 01 - 30021 N\^\i mes, France.\newline \indent I3M (Institut de Math\'ematiques et Mod\'elisation de Montpellier) UMR - CNRS 5149, Universit\'e Montpellier II, Place Eug\`ene Bataillon, 34090 Montpellier, France.}
\email{jean-philippe.mandallena@unimes.fr}

\begin{document}

\maketitle

%%%%%%%%%%%%%% INTRODUCTION %%%%%%%%%%%%%%%%%%%%

\begin{quote}
\footnotesize{\sc R\'esum\'e.} Acerbi, Buttazzo et Percivale ont donn\'es une d\'efinition variationnelle de l'\'energie d'une corde non lin\'eaire sous la contrainte ``$\det\nabla u>0$'' (voir \cite{acbuper}). Dans le m\^eme esprit, nous obtenons l'\'energie d'une membrane non lin\'eaire sous la contrainte plus simple ``$\det\nabla u\not=0$''.
\end{quote}

\bigskip

\begin{quote}
\footnotesize{\sc Abstract.} Acerbi, Buttazzo and Percivale gave a variational definition of the nonlinear string energy under the constraint ``$\det\nabla u>0$'' (see \cite{acbuper}). In the same spirit, we obtain the nonlinear membrane energy under the simpler constraint ``$\det\nabla u\not=0$''\footnotemark[1].

\vskip0.5mm

\end{quote}

\section{Introduction}

\footnotetext[1]{In \cite{benbelgacem1}, Ben Belgacem announced to have obtained a variational definition of the nonlinear membrane energy under the constraint ``$\det\nabla u>0$''. To our knowledge, his statement \cite[Theorem 1]{benbelgacem1} never was proved (see Remark \ref{benbelgacemwork}).}

Consider an elastic material occupying in a reference configuration the bounded open set $\Sigma_\eps\subset\RR^3$ given by
$$
\Sigma_\eps:=\Sigma\times\left]-{\eps\over 2},{\eps\over 2}\right[,
$$
where $\eps>0$ is very small and $\Sigma\subset\RR^2$ is Lipschitz, open and bounded. A point of $\Sigma_\eps$ is denoted by $(x,x_3)$ with $x\in\Sigma$ and $x_3\in]-{\eps\over 2},{\eps\over 2}[$.
Denote by 
$
W:\MM^{3\times 3}\to[0,+\infty]
$ 
the stored-energy function supposed to be {\em continuous} and {\em coercive}, i.e., $W(F)\geq C|F|^p$ for all $F\in\MM^{3\times 3}$ and some $C>0$. In order to take into account the fact that an infinite amount of energy is required to compress a finite volume into zero volume\footnote[2]{However, we do not prevent orientation reversal.}, i.e., 
\begin{equation}\label{det}
W(F)\to+\infty\ \hbox{ as }\ \det F\to 0,
\end{equation}
where $\det F$ denotes the determinant of the $3\times 3$ matrix $F$, we assume that:
\begin{itemize}
\item[(C$_1$)] {\em $W(F)=+\infty$ if and only if $\det F=0$};
\item[(C$_2$)] {\em for every $\delta>0$, there exists $c_\delta>0$ such that for all $F\in\MM^{3\times 3}$,}
$$
\hbox{\em if }|\det F|\geq\delta\hbox{\em\ then } W(F)\leq c_\delta(1+|F|^p). 
$$
\end{itemize}
Our goal is to show that as $\eps\to 0$ the three-dimensional free energy functional $\I_\eps:W^{1,p}(\Sigma_\eps;\RR^3)\to[0,+\infty]$ (with $p>1$) defined by
\begin{equation}\label{SEF}
\I_\eps(u):={1\over\eps}\int_{\Sigma_\eps}W(\nabla u(x,x_3))dxdx_3
\end{equation}
converges in a variational sense (see Definition \ref{variationalconvergence}) to the two-dimensional free energy functional $\I_{\rm mem}:W^{1,p}(\Sigma;\RR^3)\to[0,+\infty]$ given by
\begin{equation}\label{MSEF}
\I_{\rm mem}(v):=\int_\Sigma W_{\rm mem}(\nabla v(x))dx
\end{equation}
with $W_{\rm mem}:\MM^{3\times 2}\to[0,+\infty]$. Usually, $E_{\rm mem}$ is called the nonlinear membrane energy associated with the two-dimensional elastic material with respect to the reference configuration $\Sigma$. Furthermore we wish to give a representation formula for $W_{\rm mem}$. 

Such a problem was studied by Le Dret and Raoult in \cite{ledretraoult} when $W$ is of $p$-polynomial growth, i.e., $W(F)\leq c(1+|F|^p)$ for all $F\in\MM^{3\times 3}$ and some $c>0$, so that (\ref{det}) is not satisfied. The distinguishing feature here is that $W$ is not of $p$-polynomial growth. 

An outline of the paper is as follows. The variational convergence of $E_\eps$ to $E_{\rm mem}$ as $\eps\to 0$ as well as a representation formula for $W_{\rm mem}$ are given by Corollary \ref{corollary} (see also Proposition \ref{loads}). Corollary \ref{corollary} is a consequence of Theorems \ref{first_main_result} and \ref{anzman}. As Theorem \ref{anzman} is proved in our previous article \cite{anzman2}, the main result of the paper is Theorem \ref{first_main_result}. In fact, Theorem \ref{anzman} is analogous to Theorem \ref{benbelgacem} established by Ben Belgacem in \cite{benbelgacem}. A comparison of these results is made in Sect. 2.3 (see also \cite[Remark 2.6]{anzman2}). Theorem \ref{first_main_result} is proved in Section 4: the principal ingredients being Theorem \ref{add} (stated in Sect. 2.2 and whose proof is contained in \cite{anzman2}) and Theorem \ref{basictheorem} (whose statement and proof are given in Section 3). 

For the convenience of the reader, we recall the proofs of Theorems \ref{add} and \ref{anzman} in appendix.

%%%%%%%%%%%%%% MAIN RESULTS %%%%%%%%%%%%%%%%%%%%%

\section{Results}

                                %%%%% VARIATIONAL CONVERGENCE %%%%%%%

\subsection{Variational convergence} As in \cite{acbuper}, to accomplish our asymptotic analysis, we use the notion of convergence introduced by Anzellotti, Baldo and Percivale in \cite{anzebalper} in order to deal with dimension reduction problems in mechanics. Let $\pi=\{\pi_\eps\}_\eps$ be the family of maps $\pi_\eps:W^{1,p}(\Sigma_\eps;\RR^3)\to W^{1,p}(\Sigma;\RR^3)$ defined by
$$
\pi_\eps(u):={1\over\eps}\int_{-{\eps\over 2}}^{\eps\over 2}u(\cdot,x_3)dx_3.
$$
\begin{definition}\label{variationalconvergence}
We say that $\I_\eps$ $\Gamma(\pi)$-converges to $\I_{\rm mem}$ as $\eps\to 0$, and we write $\I_{\rm mem}=\Gamma(\pi)\hbox{\rm -}\lim_{\eps\to 0}\I_\eps$, if the following two assertions hold{\rm:}
\begin{itemize}
\item[(i)] for all $v\in W^{1,p}(\Sigma;\RR^3)$ and all $\{u_\eps\}_\eps\subset W^{1,p}(\Sigma_\eps;\RR^3)$, 
$$
\hbox{if }\pi_\eps(u_\eps)\to v\hbox{ in }L^{p}(\Sigma;\RR^3)\hbox{ then }\I_{\rm mem}(v)\leq\liminf_{\eps\to 0}\I_\eps(u_\eps);
$$
\item[(ii)] for all $v\in W^{1,p}(\Sigma;\RR^3)$, there exists $\{u_\eps\}_\eps\subset W^{1,p}(\Sigma_\eps;\RR^3)$ such that{\rm:} 
$$
\pi_\eps(u_\eps)\to v\hbox{ in }L^{p}(\Sigma;\RR^3)\hbox{ and }\I_{\rm mem}(v)\geq\limsup_{\eps\to 0}\I_\eps(u_\eps).
$$
\end{itemize}
\end{definition}
 In fact, Definition \ref{variationalconvergence} is a variant of De Giorgi's $\Gamma$-convergence. This is made clear by Lemma \ref{link}. Consider $\mathcal{\I}_\eps:W^{1,p}(\Sigma;\RR^3)\to[0,+\infty]$ defined by
$$
\mathcal{\I}_\eps(v):=\inf\Big\{\I_\eps(u):\pi_\eps(u)=v\Big\}.
$$
\begin{definition}\label{gammavariationalconvergence}
We say that $\mathcal{\I}_\eps$ $\Gamma$-converges to $\I_{\rm mem}$ as $\eps\to 0$, and we write $\I_{\rm mem}=\Gamma\hbox{\rm -}\lim_{\eps\to 0}\mathcal{\I}_\eps$, if for every $v\in W^{1,p}(\Sigma;\RR^3)$,
$$
\left(\Gamma\hbox{-}\liminf_{\eps\to 0}\mathcal{\I}_\eps\right)(v)=\left(\Gamma\hbox{-}\limsup_{\eps\to 0}\mathcal{\I}_\eps\right)(v)=E_{\rm mem}(v),
$$
where
$
\big(\Gamma\hbox{-}\liminf_{\eps\to 0}\mathcal{\I}_\eps\big)(v):=\inf\big\{\liminf_{\eps\to 0}\mathcal{\I}_\eps(v_\eps):v_\eps\to v\hbox{ in }L^{p}(\Sigma;\RR^3)\big\}
$
and
$
\big(\Gamma\hbox{-}\limsup_{\eps\to 0}\mathcal{\I}_\eps\big)(v):=\inf\big\{\limsup_{\eps\to 0}\mathcal{\I}_\eps(v_\eps):v_\eps\to v\hbox{ in }L^{p}(\Sigma;\RR^3)\big\}.
$
\end{definition}
For a deeper discussion of the $\Gamma$-convergence theory we refer to the book \cite{dalmaso}. Definition \ref{gammavariationalconvergence} is equivalent to assertions (i) and (ii) in Definition \ref{variationalconvergence} with ``$\pi(u_\eps)\to v$'' replaced by ``$v_\eps\to v$''. It is then obvious that
\begin{lemma}\label{link}
$\I_{\rm mem}=\Gamma(\pi)\hbox{\rm -}\lim_{\eps\to 0}\I_\eps$ if and only if $\I_{\rm mem}=\Gamma\hbox{\rm -}\lim_{\eps\to 0}\mathcal{\I}_\eps$.
\end{lemma}
As in \cite{acbuper}, suppose that the exterior loads derive from a potential $\Psi:\overline{\Sigma}_1\times\RR^3\to\RR$ given by 
$
\Psi((x,x_3),\zeta):=\langle \psi(x,x_3),\zeta\rangle+|\zeta|^p,
$ 
where $\psi:\overline{\Sigma}_1\to\RR^3$ is continuous and $\langle\cdot,\cdot\rangle$ denotes the scalar product in $\RR^3$, and define $L_\eps:W^{1,p}(\Sigma_\eps;\RR^3)\to\RR$ and $L_{\rm mem}:W^{1,p}(\Sigma;\RR^3)\to\RR$ by
$$
L_\eps(u):={1\over\eps}\int_{\Sigma_\eps}\Psi((x,x_3),u(x,x_3))dxdx_3\hbox{ and }L_{\rm mem}(v):=\int_{\Sigma}\Psi((x,0),v(x))dx.
$$
Then, using similar arguments to those in \cite[proof of Proposition 3.1 p. 141 and proof of Theorem 2.1 p. 145]{acbuper}, we obtain
\begin{proposition}\label{loads} 
Assume that $\I_\eps$ in {\rm(\ref{SEF})} $\Gamma(\pi)$-converge to $E_{\rm mem}$ in {\rm(\ref{MSEF})} as $\eps\to 0$, and consider $\{u_\eps\}_\eps\subset W^{1,p}(\Sigma_\eps;\RR^3)$ such that
$$
\I_\eps(u_\eps)+L_\eps(u_\eps)-\inf\Big\{\I_\eps(u)+L_\eps(u):u\in W^{1,p}(\Sigma_\eps;\RR^3)\Big\}\to 0\hbox{ as }\eps\to0.
$$
Then, $\{\pi_\eps(u_\eps)\}_\eps$ is weakly relatively compact in $W^{1,p}(\Sigma;\RR^3)$ and each of its cluster points $\bar v$ satisfies 
$$
E_{\rm mem}(\bar v)+L_{\rm mem}(\bar v)=\min\Big\{E_{\rm mem}(v)+L_{\rm mem}(v):v\in W^{1,p}(\Sigma;\RR^3)\Big\}.
$$
\end{proposition}

                                    %%%%%%% TRANSITION %%%%%%% 

The method used in this paper for passing from (\ref{SEF}) to (\ref{MSEF}) was initiated by Anza Hafsa in \cite{anza1,anza} (see also Mandallena \cite{mandallena1,mandallena2} and Anza Hafsa-Mandallena \cite{anzman1,anzman3} for the relaxation case). It first consists of studying the $\Gamma$-convergence of $\mathcal{\I}_\eps$ as $\eps\to 0$ (see Sect. 2.2), and then establishing an integral representation for the corresponding $\Gamma$-limit (see Sect. 2.3). 

                                    %%%%%% MAIN RESULT %%%%%%

\subsection{\boldmath$\Gamma$\unboldmath-convergence of \boldmath$\mathcal{\I}_\eps$\unboldmath\ as \boldmath$\eps\to 0$\unboldmath} From now on, given a bounded open set $D\subset\RR^2$ with $|\partial D|=0$, we denote by $\Aff(D;\RR^3)$ the space of all continuous piecewise affine functions from $D$ to $\RR^3$, i.e., {\em $v\in\Aff(D;\RR^3)$ if and only if $v$ is continuous and there exists a finite family $(D_i)_{i\in I}$ of open disjoint subsets of $D$ such that $|\partial D_i|=0$ for all $i\in I$, $|D\setminus \cup_{i\in I} D_i|=0$ and for every $i\in I$, $\nabla v(x)=\xi_i$ in $D_i$ with $\xi_i\in\MM^{3\times 2}$} (where $|\cdot|$ denotes the Lebesgue measure in $\RR^2$).  
\begin{remark}\label{EkelandTemamDensity}
From Ekeland-Temam \cite{ekeland}, we know that $\Aff^{ET}(D;\RR^3)$ is strongly dense in $W^{1,p}(D;\RR^3)$, where $\Aff^{ET}(D;\RR^3)$ is defined as follows: {\em $v\in\Aff^{ET}(D;\RR^3)$ if and only if $v$ is continuous and there exists a finite family $(D_i)_{i\in I}$ of open disjoint subsets of $D$ such that $|\partial D_i|=0$ for all $i\in I$, $|D\setminus \cup_{i\in I} D_i|=0$ and for every $i\in I$, the restriction of $v$ to $D_i$ is affine}. As $\Aff^{ET}(D;\RR^3)\subset\Aff(D;\RR^3)\subset W^{1,p}(D;\RR^3)$, it is clear that $\Aff(D;\RR^3)$ is also strongly dense in $W^{1,p}(D;\RR^3)$. (Note that the fact of considering $\Aff(D;\RR^3)$ instead  of $\Aff^{ET}(D;\RR^3)$ plays an important role in our analysis, see Remarks \ref{CPAFremark} and \ref{CPAFremark2}.)
\end{remark}
Let $\J:W^{1,p}(\Sigma;\RR^3)\to[0,+\infty]$ be defined by
$$
\J(v):=\left\{
\begin{array}{cl}
\displaystyle\int_\Sigma W_0(\nabla v(x))dx&\hbox{if }v\in\Aff(\Sigma;\RR^3)\\
+\infty&\hbox{otherwise,}
\end{array}
\right.
$$
where, as in \cite{ledretraoult}, $W_0:\MM^{3\times 2}\to[0,+\infty]$ is given by
$$
W_0(\xi):=\inf_{\zeta\in\RR^3}W(\xi\mid\zeta)
$$
with $(\xi\mid\zeta)$ denoting the element of $\MM^{3\times 3}$ corresponding to $(\xi,\zeta)\in\MM^{3\times 2}\times\RR^3$. (As $W$ is coercive, it is easy to see that {\em$W_0$ is coercive}, i.e., $W_0(\xi)\geq C|\xi|^p$ for all $\xi\in\MM^{3\times 2}$ and some $C>0$.) Note that conditions (C$_1$) and (C$_2$) imply $W_0$ is not of $p$-polynomial growth. In fact, we have
\begin{lemma}\label{propertiesofW_0}
Denote by $\xi_1\land\xi_2$ the cross product of vectors $\xi_1,\xi_2\in\RR^3$.
\begin{itemize}
\item[(i)] If {\rm(C$_1$)} holds then
\begin{itemize}
\item[($\overline{\rm C}_1$)] $W_0(\xi_1\mid\xi_2)=+\infty$ if and only if $\xi_1\land\xi_2=0$.
\end{itemize}
\item[(ii)] If {\rm(C$_2$)} holds then 
\begin{itemize}
\item[($\overline{\rm C}_2$)] for all $\delta>0$, there exists $c_\delta>0$ such that for all $\xi=(\xi_1\mid\xi_2)\in\MM^{3\times 2}$,
$$
\hbox{if }|\xi_1\land\xi_2|\geq\delta\hbox{ then } W_0(\xi)\leq c_\delta(1+|\xi|^p). 
$$
\end{itemize}
\end{itemize}
\end{lemma}
\begin{proof} 
(i) Given $\xi=(\xi_1\mid\xi_2)$, if $W_0(\xi_1\mid\xi_2)<+\infty$ (resp. $W_0(\xi_1\mid\xi_2)=+\infty$) then $W(\xi\mid\zeta)<+\infty$ (resp. $W(\xi\mid\zeta)=+\infty$) for some $\zeta\in\RR^3$ (resp. for all $\zeta\in\RR^3$), and so $\xi_1\land\xi_2\not=0$ (resp. $\xi_1\land\xi_2=0$) by (C$_1$).

(ii) Let $\delta>0$ and let $\xi=(\xi_1\mid\xi_2)$ be such that $|\xi_1\land\xi_2|\geq\delta$. Setting $\zeta:={\xi_1\land\xi_2\over|\xi_1\land\xi_2|}$, we have $\det (\xi\mid\zeta)\geq\delta$, and   using (C$_2$) we can assert that there exists $c_\delta>0$, which does not depend on $\xi$, such that $W_0(\xi)\leq c_\delta(1+|\xi|^p)$.
\end{proof}

Assume furthermore that 
\begin{itemize}
\item[(C$_3$)] {\em $W(\xi\mid\zeta)=W(\xi\mid-\zeta)$ for all $\xi\in\MM^{3\times 2}$ and all $\zeta\in\RR^3$.}
\end{itemize}
The main result of the paper is the following.
\begin{theorem}\label{first_main_result}
Under {\rm(C$_1$)}, {\rm(C$_2$)} and  {\rm(C$_3$)}, we have $\Gamma\hbox{\rm -}\lim_{\eps\to 0}\mathcal{\I}_\eps=\overline{\J}$ with  $\overline{\J}:W^{1,p}(\Sigma;\RR^3)\to[0,+\infty]$ given by
$$
\overline{\J}(v):=\inf\left\{\liminf_{n\to+\infty}\J(v_n):W^{1,p}(\Sigma;\RR^3)\ni v_n\to v\hbox{ in }L^{p}(\Sigma;\RR^3)\right\}.
$$
\end{theorem}
 The proof of Theorem \ref{first_main_result} is established in Section 4. It uses Theorem \ref{basictheorem} (see Section 3) and Theorem \ref{add}  whose proof is contained in \cite{anzman2}.
\begin{theorem}\label{add}
If {\rm($\overline{\rm C}_2$)} holds then $\overline{\J}=\JJ$ with $\JJ:W^{1,p}(\Sigma;\RR^3)\to[0,+\infty]$ given by
$$
\JJ(v)=\inf\left\{\liminf_{n\to+\infty}\int_\Sigma W_0(\nabla v_n(x))dx:W^{1,p}(\Sigma;\RR^3)\ni v_n\to v\hbox{ in }L^{p}(\Sigma;\RR^3)\right\}.
$$
\end{theorem}
\begin{remark}
Theorem \ref{first_main_result} can be applied when $W:\MM^{3\times 3}\to[0,+\infty]$ is given by
$$
W(F):=h(|\det F|)+|F|^p,
$$
where $h:[0,+\infty[\to[0,+\infty]$ is a continuous functions such that:
\begin{itemize}
\item[-] {$h(t)=+\infty$ if and only if $t=0$};
\item[-] {for every $\delta>0$, there exists $r_\delta>0$ such that $h(t)\leq r_\delta$ for all $t\geq\delta$}.
\end{itemize}
\end{remark}

                                                                       %%%% RELAXATION %%%%%

\subsection{Integral representation of \boldmath$\overline{\J}$\unboldmath} Our framework leads us to deal with relaxation of nonconvex integral functionals which are not of $p$-polynomial growth. Such relaxation problems were studied in Ben Belgacem \cite{benbelgacem} and Anza Hafsa-Mandallena \cite{anzman2} (see also Carbone-De Arcangelis \cite{carde} for the scalar case). To state the integral representation theorems obtained in these papers (see Theorems \ref{benbelgacem} and \ref{anzman}), we need the concepts of quasiconvex envelope  and rank-one convex envelope.
\begin{definition}\label{DEFofQUasiRankoNEConvexityandEnvelOPe}
Let $f:\MM^{3\times 2}\to[0,+\infty]$ be a Borel measurable function.
\begin{itemize}
\item[(i)] We say that $f$ is quasiconvex if for every $\xi\in\MM^{3\times 2}$, every bounded open set $D\subset\RR^2$ with $|\partial D|=0$ and every $\phi\in W^{1,\infty}_0(D;\RR^3)$,
$$
f(\xi)\leq{1\over|D|}\int_D f(\xi+\nabla\phi(x))dx.
$$ 
\item[(ii)] By the quasiconvex envelope of $f$, we mean the unique function (when it exists) $\mathcal{Q}f:\MM^{3\times 2}\to[0,+\infty]$ such that:
\begin{itemize}
\item[-] $\mathcal{Q}f$ is Borel measurable, quasiconvex and $\mathcal{Q}f\leq f$;
\item[-] for all $g:\MM^{3\times 2}\to[0,+\infty]$, if $g$ is Borel measurable, quasiconvex and $g\leq f$, then $g\leq\mathcal{Q}f$.
\end{itemize} 
(Usually, for simplicity, we say that $\mathcal{Q}f$ is the greatest quasiconvex function which less than or equal to $f$.) 
\item[(iii)] We say that $f$ is rank one convex if for every $\alpha\in]0,1[$ and every $\xi,\xi^\prime\in\MM^{3\times 2}$ with rank($\xi-\xi^\prime$)=1,
$$
f(\alpha\xi+(1-\alpha)\xi^\prime)\leq \alpha f(\xi)+(1-\alpha)f(\xi^\prime).
$$
\item[(iv)] By the rank one convex envelope of $f$, that we denote by $\mathcal{R}f$, we mean the greatest rank one convex function which less than or equal to $f$.
\end{itemize}
\end{definition}
\begin{remark}\label{rankoneconvexity}
It is well known that if $f$ is quasiconvex and {\em continuous} then $f$ is rank one convex. This is false for a general Borel measurable $f$ (see \cite[Example 3.5]{ballmurat}).
\end{remark}

                                                  %%%%%% Ben Belgacem's relaxation theorem %%%%%

\subsubsection{Ben Belgacem's theorem} In \cite[Section 5.1]{benbelgacem} Ben Belgacem asserts that if $W_0$ satisfies {\rm($\overline{\rm C}_2$)} then $\mathcal{R}W_0$ is of $p$-polynomial growth, so that is $\mathcal{Q}[\mathcal{R}W_0]$. (As $W_0$ is coercive, it is easy to see that $\mathcal{R}W_0$ is coercive.) Then, using his main result \cite[Theorem 3.1]{benbelgacem}, he obtains
\begin{theorem}\label{benbelgacem}
If {\rm($\overline{\rm C}_2$)} holds then for every $v\in W^{1,p}(\Sigma;\RR^3)$,
$$
\overline{\J}(v)=\int_\Sigma \mathcal{Q}[\mathcal{R}W_0](\nabla v(x))dx.
$$
\end{theorem}
 According to Remark \ref{rankoneconvexity}, in Theorem \ref{benbelgacem} we cannot know if $\mathcal{Q}[\mathcal{R}W_0]=\mathcal{Q}W_0$. In fact, under {\rm($\overline{\rm C}_2$)} the latter equality holds (see Remark \ref{anzmanremak}).
\begin{remark}\label{benbelgacemwork}
In \cite[Theorem 1]{benbelgacem1} Ben Belgacem announced to have established the $\Gamma(\pi)$-convergence of $E_\eps$ to $E_{\rm mem}$ as $\eps\to 0$ under the   two (more physical) conditions:
\begin{itemize}
\item[($\hat{\rm C}_1$)] {\em $W(F)=+\infty$ if and only if $\det F\leq 0$};
\item[($\hat{\rm C}_2$)] {\em for every $\delta>0$, there exists $c_\delta>0$ such that for all $F\in\MM^{3\times 3}$,}
$$
\hbox{\em if }\det F\geq\delta\hbox{\em\ then } W(F)\leq c_\delta(1+|F|^p). 
$$
\end{itemize}
In \cite{benbelgacem}, which is the paper corresponding to the note \cite{benbelgacem1}, the statement \cite[Theorem 1]{benbelgacem1} is not proved. To our knowledge, under ($\hat{\rm C}_1$) and ($\hat{\rm C}_2$) the problem of passing from (\ref{SEF}) to (\ref{MSEF}) by using $\Gamma(\pi)$-convergence is still open.
\end{remark}
                
           %%%%%%% An alternative relaxation theorem %%%%%%

\subsubsection{An alternative theorem} Define $\Z W_0:\MM^{3\times 2}\to[0,+\infty]$ by 
\begin{equation}\label{DefinitionofZW0}
\Z W_0(\xi):=\inf\left\{\int_Y W_0\big(\xi+\nabla\phi(y)\big)dy:\phi\in \Aff_0(Y;\RR^3)\right\}.
\end{equation}
with $Y:=]0,1[^2$ and $\Aff_0(D;\RR^3):=\{\phi\in\Aff(Y;\RR^3):\phi=0\hbox{ on }Y\}$. (As $W_0$ is coercive, it is easy to see that {\em$\Z W_0$ is coercive}.) In \cite{anzman2}, under  {\rm($\overline{\rm C}_2$)}, we prove that  $\Z W_0$ is of $p$-polynomial growth and continuous (see Propositions \ref{theorAA2} and \ref{propAA1}(iii)), and that $\Z W_0$ is the quasiconvex envelope of $W_0$, i.e., $\Z W_0=\mathcal{Q}W_0$ (see Propositions \ref{corA4}). Theorem \ref{anzman}  is contained in \cite{anzman2} (for the convenience of the reader, we give the proof in appendix).
\begin{theorem}\label{anzman}
If {\rm($\overline{\rm C}_2$)} holds then for every $v\in W^{1,p}(\Sigma;\RR^3)$,
$$
\overline{\J}(v)=\int_\Sigma \mathcal{Q}W_0(\nabla v(x))dx.
$$
\end{theorem}
\begin{remark}\label{anzmanremak}
If ($\overline{\rm C}_2$) holds then $\mathcal{Q}[\mathcal{R}W_0]=\mathcal{Q}W_0$. Indeed, by Proposition \ref{theorAA2}, $\Z W_0(\xi)\leq c(1+|\xi|^p)$ for all $\xi\in\MM^{3\times 2}$ and some $c>0$. Then $\Z W_0$ is finite, and so  $\Z W_0$ is continuous by Proposition \ref{propAA1}(iii). It follows that $\Z W_0=\mathcal{Q}W_0$ (see the proof of Proposition  \ref{corA4}). Thus $\mathcal{Q}W_0$ is continuous, hence  $\mathcal{Q}W_0$ is rank-one convex (see Remark \ref{rankoneconvexity}), and the result follows.
\end{remark}

\subsection{\boldmath$\Gamma(\pi)$\unboldmath-convergence of \boldmath$\I_\eps$\unboldmath\ to \boldmath$\I_{\rm mem}$\unboldmath\ as \boldmath$\eps\to 0$\unboldmath} According to Lemmas \ref{link} and \ref{propertiesofW_0}(ii), a direct consequence of Theorems \ref{first_main_result} and \ref{anzman} is the following.
\begin{corollary}\label{corollary}
Let assumptions {\rm(C$_1$)}, {\rm(C$_2$)} and {\rm(C$_3$)}  hold. Then as $\eps\to 0$, $\I_\eps$ in {\rm(\ref{SEF})} $\Gamma(\pi)$-converge to $\I_{\rm mem}$ in {\rm(\ref{MSEF})} with $W_{\rm mem}=\mathcal{Q}W_0$. 
\end{corollary}

%%%%% A representation formula for \boldmath$\J$\unboldmath %%%%%%%%%%%

\section{Representation of $\J$} 

The goal of this section is to show Theorem \ref{basictheorem}. To this end, we begin by proving three lemmas. From now on, we set
$$
\Aff_*(\Sigma;\RR^3):=\Big\{v\in\Aff(\Sigma;\RR^3):\partial_1v(x)\land\partial_2 v(x)\not=0\hbox{ a.e. in }\Sigma\Big\},
$$
where $\partial_1v(x)$ (resp. $\partial_2 v(x)$) denotes the partial derivative of $v$ at $x=(x_1,x_2)$ with respect to $x_1$ (resp. $x_2$). By definition, to every $v\in\Aff_*(\Sigma;\RR^3)$ there corresponds a finite family $(V_i)_{i\in I}$ of open disjoint subsets of $\Sigma$ such that:
\begin{itemize}
\item[-] $|\partial V_i|=0$ for all $i\in I$;
\item[-] $|\Sigma\setminus\cup_{i\in I}V_i|=0$; 
\item[-] for every $i\in I$, $\nabla v(x)=\xi_{i}$ in $V_i$ with $\xi_{i}=(\xi_{i,1}\mid\xi_{i,2})\in\MM^{3\times 2}$; 
\item[-] $\xi_{i,1}\land\xi_{i,2}\not=0$ for all $i\in I$. 
\end{itemize}
\begin{lemma}\label{Lemma0}
If {\rm (C$_1$)} holds then $\dom\J=\Aff_*(\Sigma;\RR^3)$, where $\dom\J$ is the effective domain of $\J$.
\end{lemma}
\begin{proof}It is a direct consequence of Lemma \ref{propertiesofW_0}(i).
\end{proof}

Given $v\in\Aff_*(\Sigma;\RR^3)$, for every $i\in I$ and every integer $j\geq 1$, we consider the subsets $U^-_{i,j}$ and $U^+_{i,j}$ of $\RR^3$ given by
$$
U_{i,j}^-:=\left\{\zeta\in\RR^3:\det(\xi_{i}\mid\zeta)\leq-{1\over j}\right\}\hbox{ and } U_{i,j}^+:=\left\{\zeta\in\RR^3:\det(\xi_{i}\mid\zeta)\geq {1\over j}\right\}.
$$
Here are some elementary properties of these sets:
\begin{itemize}
\item[(P$_1$)] {both $U_{i,j}^-$ and $U_{i,j}^+$ are nonempty convex subsets of $\RR^3$};
\item[(P$_2$)] $U_{i,j}^-\cup U_{i,j}^+=\big\{\zeta\in\RR^3:|\det(\xi_i\mid\zeta)|\geq {1\over j}\big\}$;
\item[(P$_3$)] $U_{i,1}^-\subset U_{i,2}^-\subset U_{i,3}^-\subset\cdots\subset\union_{j\geq 1}U_{i,j}^-=\left\{\zeta\in\RR^3:\det\big(\xi_i\mid\zeta\big)< 0\right\}$;
\item[(P$_4$)] $U_{i,1}^+\subset U_{i,2}^+\subset U_{i,3}^+\subset\cdots\subset\union_{j\geq 1}U_{i,j}^+=\left\{\zeta\in\RR^3:\det\big(\xi_i\mid\zeta\big)> 0\right\}$.
\end{itemize} 
\begin{lemma}\label{Lemma1}
Given $v\in\Aff_*(\Sigma;\RR^3)$, there exist $j_v\geq 1$ and two subsets $I^-$ and $I^+$ of $I$, with $I^-\cup I^+=I$ and $I^-\cap I^+=\emptyset$, such that  for all $j\geq j_v$,
$$
\left(\Cap_{i\in I^-}U_{i,j}^-\right)\cap\left(\Cap_{i\in I^+}U_{i,j}^+\right)\not=\emptyset.
$$
\end{lemma}
\begin{proof}
For every $i\in I$, define the hyperplane $H_i$ of $\RR^3$ by
$
H_i:=\{\zeta\in\RR^3:\det(\xi_i\mid \zeta)=0\}.
$
It is obvious that $\cup_{i\in I}H_i\not=\RR^3$, and so there exists $\zeta\in\RR^3$ such that
$
\det(\xi_i\mid \zeta)\not=0
$
for all $i\in I$. Taking (P$_2$) into account, we deduce the existence of an integer $j_v\geq 1$ for which $\zeta\in\cap_{i\in I}\big(U_{i,j_v}^-\cup U_{i,j_v}^+\big)$. Hence, there are two subsets $I^-$ and $I^+$ of $I$, with $I^-\cup I^+=I$ and $I^-\cap I^+=\emptyset$, such that 
$
(\Cap_{i\in I^-}U_{i,j_v}^-)\cap(\Cap_{i\in I^+}U_{i,j_v}^+)\not=\emptyset,
$
and the lemma follows by using (P$_3$) and (P$_4$).
\end{proof}

Setting
$
V:=\Cup_{i\in I}V_i,
$ 
 for every $j\geq j_v$, with $j_v$ given by Lemma \ref{Lemma1}, we define $\Gamma^j_v:\overline{\Sigma}\dto\RR^3$ by
$$
\Gamma_v^j(x):=\left\{
\begin{array}{cl}
U_{i,j}^-&\hbox{ if }x\in V_i\hbox{ with }i\in I^-\\
U_{i,j}^+&\hbox{ if }x\in V_i\hbox{ with }i\in I^+\\
\left(\Cap\limits_{i\in I^-}U_{i,j}^-\right)\cap\left(\Cap\limits_{i\in I^+}U_{i,j}^+\right)&\hbox{ if }x\in \overline{\Sigma}\setminus V.
\end{array}
\right.
$$
(It is clear that for every $x\in\overline{\Sigma}$, $\Gamma_v^j(x)$ is a nonempty convex closed subset of $\RR^3$.) In the sequel, given $\Gamma:\overline{\Sigma}\dto\RR^3$ we set
$$
C(\overline{\Sigma};\Gamma):=\Big\{\phi\in C(\overline{\Sigma};\RR^3):\phi(x)\in\Gamma(x)\hbox{ a.e. in }\overline{\Sigma}\Big\},
$$ 
where $C(\overline{\Sigma};\RR^3)$ denotes the space of all continuous functions from $\overline{\Sigma}$ to $\RR^3$. 
\begin{lemma}\label{Lemma3}
Given $v\in\Aff_*(\Sigma;\RR^3)$ and $j\geq j_v$, if {\rm(C$_2$)}  holds, then
$$
\inf_{\phi\in C(\overline{\Sigma};\Gamma_v^j)}\int_{\Sigma}W(\nabla v(x)\mid\phi(x))dx=\int_\Sigma \inf_{\zeta\in\Gamma_v^j(x)}W(\nabla v(x)\mid\zeta)dx.
$$
\end{lemma} 
\begin{proof}
It is obvious that
$$
\inf_{\phi\in C(\overline{\Sigma};\Gamma_v^j)}\int_{\Sigma}W(\nabla v(x)\mid\phi(x))dx\geq\int_\Sigma \inf_{\zeta\in\Gamma_v^j(x)}W(\nabla v(x)\mid\zeta)dx.
$$
Prove then the converse inequality. By Lemma \ref{Lemma1}, $(\Cap_{i\in I^-}U_{i,j}^-)\cap(\Cap_{i\in I^+}U_{i,j}^+)\not=\emptyset$, and so there exists $\bar\zeta\in(\Cap_{i\in I^-}U_{i,j}^-)\cap(\Cap_{i\in I^+}U_{i,j}^+)$. As every $U^-_{i,j}$ (resp. $U^+_{i,j}$) is closed, every $W(\xi_i\mid\cdot)$ is continuous and $W$ is coercive, for each $i\in I^-$ (resp. $i\in I^+$), there exists $\zeta_i\in U^-_{i,j}$ (resp. $\zeta_i\in U^+_{i,j}$) such that
\begin{equation}\label{AddING}
W(\xi_i\mid\zeta_i)=\inf_{\zeta\in U^-_{i,j}}W(\xi\mid\zeta)\quad (\hbox{resp. }W(\xi_i\mid\zeta_i)=\inf_{\zeta\in U^+_{i,j}}W(\xi\mid\zeta)).
\end{equation}
Fix any $n\geq 1$. Consider $\alpha_n:\overline{\Sigma}\to\RR$ given by $\alpha_n(x):=h(n{\rm dist}(x,\overline{\Sigma}\setminus V))$, where ${\rm dist}(x,\overline{\Sigma}\setminus V):=\inf\{|x-y|:y\in \overline{\Sigma}\setminus V\}$ and $h:[0,+\infty[\to[0,1]$ is a continuous function such that $h(0)=0$ and $h(t)=1$ for all $t\geq 1$. Define $\phi_n:\overline{\Sigma}\to\RR$ by
$$
\phi_n(x):=(1-\alpha_n(x))\bar\zeta+\alpha_n(x)\zeta_i.
$$
Clearly, $\phi_n$ is continuous and $\phi_n(x)\in\Gamma^j_v(x)$ for all $x\in\overline{\Sigma}$ since $\Gamma^j_v(x)$ is convex, and so $\phi_n\in C(\overline{\Sigma};\Gamma^j_v)$. Using (C$_2$) we deduce that $\sup_{n\geq 1}W(\nabla v(\cdot)\mid\phi_n(\cdot))\in L^1(\Sigma)$. Recalling that $W$ is continuous and taking (\ref{AddING}) into account, it is easy to see that $\lim_{n\to+\infty}W(\nabla v(x)\mid\phi_n(x))=\inf_{\zeta\in\Gamma^j_v(x)}W(\nabla v(x)\mid\zeta)$ for a.e. $x\in\Sigma$. Hence
\begin{eqnarray*}
\inf_{\phi\in C(\overline{\Sigma};\Gamma_v^j)}\int_{\Sigma}W(\nabla v(x)\mid\phi(x))dx&\leq&\lim_{n\to+\infty}\int_\Sigma W(\nabla v(x)\mid\phi_n(x))dx\\
&=&\int_\Sigma\inf_{\zeta\in\Gamma^j_v(x)}W(\nabla v(x)\mid\zeta)dx
\end{eqnarray*}
by Lebesgue's dominated convergence theorem, and the proof is complete.
\end{proof}
                         %%%% Representation formula for $\J$ %%%%

For every $j\geq j_v$, we define $\Lambda^j_v:\overline{\Sigma}\dto\RR^3$ by
$$
\Lambda_v^j(x):=\left\{
\begin{array}{cl}
U_{i,j}^-\cup U_{i,j}^+&\hbox{ if }x\in V_i\\
\Gamma^j_v(x)&\hbox{ if }x\in \overline{\Sigma}\setminus V.
\end{array}
\right.
$$
Here is our (non integral) representation theorem for $\J$.
\begin{theorem}\label{basictheorem} If {\rm(C$_1$)}, {\rm(C$_2$)} and {\rm(C$_3$)}  hold, then for every $v\in\dom\J$,
\begin{equation}\label{NIRF}
\J(v)=\inf_{j\geq j_v}\inf_{\phi\in C(\overline{\Sigma};\Lambda_v^j)}\int_{\Sigma}W(\nabla v(x)\mid\phi(x))dx.
\end{equation}
\end{theorem}
\begin{proof}
By Lemma \ref{Lemma0}, $\dom\J=\Aff_*(\Sigma;\RR^3)$. Fix $v\in\Aff_*(\Sigma;\RR^3)$ and denote by $\hat\J(v)$ the right-hand side of (\ref{NIRF}). It is easy to verify that $\J(v)\leq\hat\J(v)$. We are thus reduced to prove that
\begin{equation}\label{ineQuality}
\hat\J(v)\leq\J(v).
\end{equation}
From (C$_3$) we see that for every $j\geq j_v$ and every $x\in{\Sigma}$,
\begin{equation}\label{equality1}
\inf_{\zeta\in\Gamma^j_v(x)}W(\nabla v(x)\mid\zeta)=\inf_{\zeta\in\Lambda^j_v(x)}W(\nabla v(x)\mid\zeta).
\end{equation}
Noticing that $\Gamma^j_v(x)\subset\Lambda^j_v(x)$ for all $x\in{\Sigma}$ and using Lemma \ref{Lemma3} together with (\ref{equality1}), we obtain
\begin{equation}\label{inequality1}
\hat\J(v)\leq\inf_{j\geq j_v}\int_{\Sigma}\inf_{\zeta\in\Lambda^j_v(x)}W(\nabla v(x)\mid\zeta)dx.
\end{equation}
On the other hand, $\inf_{\zeta\in\Lambda^{j_v}_v(\cdot)}W(\nabla v(\cdot)\mid\zeta)\in L^1(\Sigma)$ by (C$_2$), and from (P$_3$) and (P$_4$) we deduce that if $x\in V$ then $\Lambda^{j_v}_v(x)\subset\Lambda^{j_v+1}_v(x)\subset\cdots\subset\cup_{j\geq j_v}\Lambda^j_v(x)$ with
$
\union_{j\geq j_v}\Lambda^j_v(x)=\{\zeta\in\RR^3:\det(\nabla v(x)\mid\zeta)\not=0\}.
$
Hence $\{\inf_{\zeta\in\Lambda^j_v(\cdot)}W(\nabla v(\cdot)\mid\zeta)\}_{j\geq j_v}$ is non-increasing and for every  $x\in V$,
\begin{equation}\label{equality}
\inf_{j\geq j_v}\inf_{\zeta\in\Lambda^j_v(x)}W(\nabla v(x)\mid\zeta)=W_0(\nabla v(x)),
\end{equation}
and (\ref{ineQuality}) follows from (\ref{inequality1}) and ({\ref{equality}}) by using the monotone convergence theorem.
\end{proof}

%%%%% Proof of Theorem \ref{first_main_result} %%%%%%%%%%%%%%%%%%

\section{Proof of Theorem \ref{first_main_result}} 

In this section we prove Theorem \ref{first_main_result}. Since $\Gamma\hbox{-}\liminf_{\eps\to 0}\mathcal{\I}_\eps\leq\Gamma\hbox{-}\limsup_{\eps\to 0}\mathcal{\I}_\eps$, we only need to show that:
\begin{itemize}
\item[(a)]$\displaystyle\overline{\J}\leq\Gamma\hbox{-}\liminf_{\eps\to 0}\mathcal{\I}_\eps$;
\item[(b)]$\displaystyle\Gamma\hbox{-}\limsup_{\eps\to 0}\mathcal{\I}_\eps\leq\overline{\J}$.
\end{itemize}
In the sequel, we follow the notation used in Section 3.

\subsection{Proof of (a)} Let $v\in W^{1,p}(\Sigma;\RR^3)$ and let  $\{v_\eps\}_\eps\subset W^{1,p}(\Sigma;\RR^3)$ be such that $v_\eps\to v$ in $L^{p}(\Sigma;\RR^3)$. We have to prove that
\begin{equation}\label{main_ineq1}
\liminf_{\eps\to 0}\J_\eps(v_\eps)\geq \overline{\J}(v).
\end{equation}
Without loss of generality we can assume that $\sup_{\eps>0}\J_\eps(v_\eps)<+\infty$. To every $\eps>0$ there corresponds $u_\eps\in\pi_{\eps}^{-1}(v_\eps)$ such that 
\begin{equation}\label{main_ineq2}
\E_\eps(v_\eps)\geq E_\eps(u_\eps)-\eps.
\end{equation}
Defining $\hat u_\eps:\Sigma_1\to\RR^3$ by $\hat u_\eps(x,x_3):=u_\eps(x,\eps x_3)$ we have
\begin{equation}\label{main_ineq3}
E_\eps\big(u_\eps\big)=\int_{\Sigma_1}W\Big(\partial_{1} \hat u_\eps(x,x_3)\mid\partial_{2} \hat u_\eps(x,x_3)\mid{1\over\eps}\partial_3 \hat u_\eps(x,x_3)\Big)dxdx_3.
\end{equation}
Using the coercivity of $W$, we deduce that $\left\|{\partial_3 \hat u_\eps}\right\|_{L^p(\Sigma_1;\RR^3)}\le c\eps^{p}$ for all $\eps>0$ and some $c>0$, and so
$
\|\hat u_\eps-v_\eps\|_{L^p(\Sigma_1;\RR^3)}\leq c^\prime\eps^p
$
by Poincar\'e-Wirtinger's inequality, where $c^\prime>0$ is a constant which does not depend on $\eps$. It follows that $\hat u_\eps\to v$ in $L^p(\Sigma_1;\RR^3)$.  For $x_3\in]-{1\over 2},{1\over 2}[$, let $w_\eps^{x_3}\in W^{1,p}(\Sigma;\RR^3)$ given by $w_\eps^{x_3}(x):=\hat u_\eps(x,x_3)$. Then (up to a subsequence) $w_\eps^{x_3}\to v$ in $L^p(\Sigma;\RR^3)$ for a.e. $x_3\in ]-{1\over 2},{1\over 2}[$. Taking (\ref{main_ineq2}) and (\ref{main_ineq3}) into account and using Fatou's lemma, we obtain
$$
\liminf_{\eps\to 0}{\J}_\eps(v_\eps)\geq\int_{-{1\over 2}}^{1\over 2}\left(\liminf_{\eps\to 0}\int_\Sigma W_0(\nabla w_\eps^{x_3}(x))dx\right)dx_3,
$$
and so $\liminf_{\eps\to 0}{\J}_\eps(v_\eps)\geq\JJ(v)$, and (\ref{main_ineq1}) follows by using Theorem \ref{add}.\hfill$\square$

\subsection{Proof of (b)} By Lemma \ref{Lemma0}, $\dom\J=\Aff_*(\Sigma;\RR^3)$. As $\Gamma\hbox{-}\limsup_{\eps\to 0}\mathcal{\I}_\eps$ is lower semicontinuous with respect to the strong topology of $L^{p}(\Sigma;\RR^3)$ (see \cite[Proposition 6.8 p. 57]{dalmaso}), it is sufficient to prove that for every $v\in\Aff_*(\Sigma;\RR^3)$,
\begin{equation}\label{limsupequality}
\limsup_{\eps\to 0}\mathcal{\I}_\eps(v)\leq\mathcal{\I}(v).
\end{equation}  
Given $v\in\Aff_*(\Sigma;\RR^3)$, fix any $j\geq j_v$ (with $j_v$ given by Lemma \ref{Lemma1}) and any $n\geq 1$. Using Theorem \ref{basictheorem} we obtain the existence of $\phi\in C(\overline{\Sigma};\Lambda^j_v)$ such that
\begin{equation}\label{mediainequality}
\int_\Sigma W(\nabla v(x)\mid\phi(x))dx\leq\J(v)+{1\over n}.
\end{equation}
By Stone-Weierstrass's approximation theorem, there exists $\{\phi_k\}_{k\geq 1}\subset C^\infty(\overline{\Sigma};\RR^3)$ such that
\begin{equation}\label{uniformityconvergence}
\phi_k\to\phi\hbox{ uniformly as }k\to+\infty.
\end{equation}
We claim that:
\begin{itemize}
\item[(c$_1$)] $|\det(\nabla v(x)\mid\phi_k(x))|\geq{1\over 2j}$ for all $x\in V$, all $k\geq k_v$ and some $k_v\geq 1$;
\item[(c$_2$)] $\lim\limits_{k\to+\infty}\int_\Sigma W(\nabla v(x)\mid\phi_k(x))dx=\int_\Sigma W(\nabla v(x)\mid\phi(x))dx$.
\end{itemize}
Indeed, setting $\mu_v:=\max_{i\in I}|\xi_{i,1}\land\xi_{i,2}|$ ($\mu_v>0$) and using (\ref{uniformityconvergence}), we deduce that there exists $k_v\geq 1$ such that for every $k\geq k_v$,
\begin{equation}\label{supequality}
\sup_{x\in\overline{\Sigma}}|\phi_k(x)-\phi(x)|<{1\over 2j\mu_v}.
\end{equation}
Let $x\in V_i$ with $i\in I$, and let $k\geq k_v$. As $\phi\in C(\overline{\Sigma};\Lambda^j_v)$ we have
\begin{equation}\label{supequality1}
|\det(\xi_i\mid\phi_k(x))|\geq{1\over j}-|\det(\xi_i\mid\phi_k(x)-\phi(x))|.
\end{equation}
Noticing that $|\det(\xi_i\mid\phi_k(x)-\phi(x))|\leq|\xi_{i,1}\land\xi_{i,2}||\phi_k(x)-\phi(x)|$, from (\ref{supequality}) and (\ref{supequality1}) we deduce that 
$
|\det(\xi_i\mid\phi_k(x))|\geq{1\over 2j},
$
and (c$_1$) is proved. Combining (c$_1$) with (C$_2$) we see that 
$
\sup_{k\geq k_v}W(\nabla v(\cdot)\mid\phi_k(\cdot))\in L^1(\Sigma).
$ 
As $W$ is continuous we have
$
\lim_{k\to+\infty}W(\nabla v(x)\mid\phi_k(x))=W(\nabla v(x)\mid\phi(x))
$
for all $x\in V$, and (c$_2$) follows by using Lebesgue's dominated convergence theorem, which completes the claim.

Fix any $k\geq k_v$ and define $\theta:]-{1\over 2},{1\over 2}[\to\RR$ by 
$
\theta(x_3):=\min_{i\in I}\inf_{x\in \overline{V}_i}|\det(\xi_i+x_3\nabla\phi_k(x)\mid\phi_k(x))|.
$
Clearly $\theta$ is continuous. By (c$_1$) we have $\theta(0)\geq{1\over 2j}$, and so there exists $\eta_v\in]0,{1\over 2}[$ such that $\theta(x_3)\geq{1\over 4j}$ for all $x_3\in]-\eta_v,\eta_v[$. Let $u_k:\Sigma_1\to\RR$ be given by
$
u_k(x,x_3):=v(x)+x_3\phi_k(x).
$
From the above it follows that
\begin{itemize}
 \item[(c$_3$)] $|\det\nabla u_k(x,\eps x_3)|\geq{1\over 4j}$ for all $\eps\in]0,\eta_v[$ and all $(x,x_3)\in V\times]-{1\over 2},{1\over 2}[$. 
\end{itemize}
As in the proof of (c$_2$), from (c$_3$) together with (C$_2$) and the continuity of $W$, we obtain
\begin{equation}\label{finalequality}
\lim_{\eps\to 0}\I_\eps(u_k)=\lim_{\eps\to 0}\int_{\Sigma_1} W(\nabla u_k(x,\eps x_3))dxdx_3=\int_\Sigma W(\nabla v(x)\mid\phi_k(x))dx.
\end{equation}

For every $\eps>0$ and every $k\geq k_v$, since $\pi_\eps(u_k)=v$ we have $\J_\eps(v)\leq \I_\eps(u_k)$. Using (\ref{finalequality}), (c$_2$) and (\ref{mediainequality}), we deduce that
$$
\limsup_{\eps\to 0}\J_\eps(v)\leq\J(v)+{1\over n},
$$
and (\ref{limsupequality}) follows by letting $n\to+\infty$.\hfill$\square$

%%%%%%%%%%%%% COMPLEMENT %%%%%%%%%%%%%%%%%%%%%%

\appendix

\section{Representation of $\overline{\J}$}

Theorems \ref{add}  and \ref{anzman} are contained in \cite{anzman2}. For the convenience of the reader, we give the proofs in this appendix. 

\subsection{Preliminary results} Throughout this appendix we will use Proposition \ref{propAA1} which gives three interesting  properties of $\Z W_0:\MM^{3\times 2}\to[0,+\infty]$ defined by (\ref{DefinitionofZW0}). The proof can be adapted from Fonseca \cite[Lemma 2.16, Lemma 2.20, Theorem 2.17 and Proposition 2.3]{fonseca} (the detailed verification is left to the reader).
\begin{proposition}\label{propAA1}
\begin{itemize} 
\item[(i)] For every bounded open set $D\subset\RR^2$ with $|\partial D|=0$ and every $\xi\in\MM^{3\times 2}$,
$$
\Z W_0(\xi)=\inf\left\{{1\over |D|}\int_D W_0(\xi+\nabla\phi(y))dy:\phi\in\Aff_0(D;\RR^3)\right\}.
$$
\item[(ii)] For every bounded open set $D\subset\RR^2$ with $|\partial D|=0$, every $\xi\in\MM^{3\times 2}$ and every $\phi\in\Aff_0(D;\RR^3)$,
$$
\Z W_0(\xi)\leq{1\over |D|}\int_D \Z W_0(\xi+\nabla\phi(x))dx.
$$
\item[(iii)] If $\Z W_0$ is finite then $\Z W_0$ is continuous.
\end{itemize}
\end{proposition}
\begin{remark}\label{CPAFremark}
In \cite{fonseca}, Fonseca proved that $Z W_0:\MM^{3\times 2}\to[0,+\infty]$ defined by
$$
Z W_0(\xi):=\inf\left\{\int_Y W_0(\xi+\nabla\phi(y))dy:\phi\in W^{1,\infty}_0(Y;\RR^3)\right\},
$$
where $W^{1,\infty}_0(Y;\RR^3):=\{\phi\in W^{1,\infty}(Y;\RR^3):\phi=0\hbox{ on }Y\}$, satisfies the  three properties:
\begin{itemize} 
\item[(j)] (\cite[Lemma 2.16]{fonseca}) for every bounded open set $D\subset\RR^2$ with $|\partial D|=0$ and every $\xi\in\MM^{3\times 2}$,
$$
Z W_0(\xi)=\inf\left\{{1\over |D|}\int_D W_0(\xi+\nabla\phi(y))dy:\phi\in W^{1,\infty}_0(D;\RR^3)\right\};
$$
\item[(jj)] (\cite[Lemma 2.20]{fonseca}) for every bounded open set $D\subset\RR^2$ with $|\partial D|=0$, every $\xi\in\MM^{3\times 2}$ and every $\phi\in\Aff^{ET}_0(D;\RR^3):=\{\phi\in\Aff^{ET}(D;\RR^3):\phi=0\hbox{ on }D\}$ (with $\Aff^{ET}(D;\RR^3)$ defined in Remark \ref{EkelandTemamDensity}),
$$
Z W_0(\xi)\leq{1\over |D|}\int_D Z W_0(\xi+\nabla\phi(x))dx;
$$
\item[(jjj)] (\cite[Theorem 2.17 and Proposition 2.3]{fonseca}) if $Z W_0$ is finite then $Z W_0$ is continuous.
\end{itemize}
The proof of (j) requires Vitali's covering theorem. Thus, by an examination of the details, we see that Proposition \ref{propAA1}(i) can be established  by following the same method as in \cite{fonseca} if $\Aff_0(D;\RR^3)$, where $D\subset\RR^2$ is a bounded open set such that $|\partial D|=0$, satisfies the  ``stability" condition: 
\begin{itemize}
\item[(S)] {\em for every $\phi\in\Aff_0(D;\RR^3)$, every bounded open set $E\subset\RR^2$ with $|\partial E|=0$ and every finite or countable family $(a_i+\alpha_i E)_{i\in I}$ of disjoint subsets of $D$ with $a_i\in\RR^3$, $\alpha_i>0$ and  $|D\setminus \cup_{i\in I}(a_i+\alpha_iE)|=0$, the function $v:D\to\RR^3$ defined by 
$$
v(x)=\alpha_i\phi\left({x-a_i\over\alpha_i}\right)\hbox{ if }x\in a_i+\alpha_i E
$$
belongs to $\Aff_0(D;\RR^3)$.}
\end{itemize} 
In fact, $\Aff_0(D;\RR^3)$ has this property, and so Proposition \ref{propAA1}(i) holds. Ben Belgacem was the first to point out the importance of considering a ``good" space of continuous piecewise affine functions. In a similar context (see \cite{benbelgacem}), he introduced the space $\Aff^V(D;\RR^3)$ of Vitali continuous piecewise affine functions as follows:  {\em $\phi\in\Aff^V(D;\RR^3)$ if and only if $\phi$ is continuous and there exists a finite or countable family $(O_i)_{i\in I}$ of disjoint open subsets of $D$ such that $|\partial O_i|=0$ for all $i\in I$, $|D\setminus\cup_{i\in I}O_i|=0$, and $\phi(x)=\xi_i\cdot x+a_i$ if $x\in O_i$, where $a_i\in\RR^3$, $\xi_i\in\MM^{3\times 2}$ and ${\rm Card}\{\xi_i:i\in I\}$ is finite} (setting $D_i:=\{x\in\cup_{i\in I}O_i:\nabla\phi(x)=\xi_i\}$ for all $i\in I$, we see that ${\rm Card}\{D_i:i\in I\}$ is finite, and so $\Aff^V(D;\RR^3)\subset\Aff(D;\RR^3)$). Clearly, $\Aff_0^V(D;\RR^3):=\{\phi\in\Aff^V(D;\RR^3):\phi=0\hbox{ on }D\}$ satisfies (S). Ben Belgacem then proved (j) replacing ``$W^{1,\infty}_0$" by ``$\Aff^V_0$". As noticed by him, since $\Aff^{ET}_0(D;\RR^3)$ does not satisfy (S), if we consider  ``$\Aff^{ET}_0$" instead of ``$W^{1,\infty}_0$", (j) seems to be false. Moreover, as the proofs of (jj) and (jjj)  need (j), if we replace ``$W^{1,\infty}_0$" by ``$\Aff_0^{ET}$", we are no longer sure that (jj) and (jjj) are true. However, Ben Belgacem also showed that these properties remain valid if we consider ``$\Aff_0^V$" instead of ``$W^{1,\infty}_0$" and ``$\Aff_0^{ET}$". As in \cite{benbelgacem}, by carefully checking, we see that the proofs given in \cite{fonseca} can be adapted to establish Proposition \ref{propAA1}(ii) and (iii).
\end{remark}
To prove Theorems \ref{add} and \ref{anzman} we will need the following proposition.
\begin{proposition}\label{theorAA2}
If {\rm($\overline{\rm C}_2$)} holds then $\Z W_0(\xi)\leq c(1+|\xi|^p)$ for all $\xi\in\MM^{3\times 2}$ and some $c>0$.
\end{proposition}
 To show Proposition \ref{theorAA2} we need the following lemma.

                  %%%%% LEMMA %%%%%

\begin{lemma}\label{lemmaA3}
If {\rm($\overline{\rm C}_2$)} holds then for every $\delta>0$, there exists $r_\delta>0$ such that for every $\xi=(\xi_1\mid\xi_2)\in\MM^{3\times 2}$,
$$
\hbox{ if }\min\{|\xi_1+\xi_2|,|\xi_1-\xi_2|\}\geq\delta\hbox{ then }\Z W_0(\xi)\leq r_\delta(1+|\xi|^p).
$$
\end{lemma}
\begin{proof}
Let $\delta>0$ and $\xi=(\xi_1\mid\xi_2)\in\MM^{3\times 2}$ be such that $\min\{|\xi_1+\xi_2|,|\xi_1-\xi_2|\}\geq\delta$. Then, one the three possibilities holds:
\begin{itemize}
\item[(i)] $|\xi_1\land\xi_2|\not=0$;
\item[(ii)] $|\xi_1\land\xi_2|=0$ with $\xi_1\not=0$;
\item[(iii)] $|\xi_1\land\xi_2|=0$ with $\xi_2\not=0$.
\end{itemize}
Set 
$
D:=\{(x_1,x_2)\in\RR^2:x_1-1<x_2<x_1+1\hbox{ and }-x_1-1<x_2<1-x_1\}
$ 
and, for each $t\in\RR$, define $\varphi_t\in\Aff_0(D;\RR)$ by
$$
\varphi_t(x_1,x_2):=\left\{
\begin{array}{ll}
-tx_1+t(x_2+1)&\hbox{if }(x_1,x_2)\in\Delta_1\\
t(1-x_1)-tx_2&\hbox{if }(x_1,x_2)\in\Delta_2\\
tx_1+t(1-x_2)&\hbox{if }(x_1,x_2)\in\Delta_3\\
t(x_1+1)+tx_2&\hbox{if }(x_1,x_2)\in\Delta_4
\end{array}
\right.
$$
with
\begin{itemize} 
\item[]$
\Delta_1:=\{(x_1,x_2)\in D:x_1\geq 0\hbox{ and } x_2\leq 0\};
$
\item[]$
\Delta_2:=\{(x_1,x_2)\in D:x_1\geq 0\hbox{ and }x_2\geq 0\};
$
\item[]$
\Delta_3:=\{(x_1,x_2)\in D:x_1\leq 0\hbox{ and }x_2\geq 0\};
$
\item[]$
\Delta_4:=\{(x_1,x_2)\in D:x_1\leq 0\hbox{ and }x_2\leq 0\}.
$
\end{itemize}
Consider $\phi\in\Aff_0(D;\RR^3)$ given by 
$$
\phi:=(\varphi_{\nu_1},\varphi_{\nu_2},\varphi_{\nu_3})\hbox{ with }\left\{\begin{array}{ll}\nu={\xi_1\land \xi_2\over|\xi_1\land \xi_2|}&\hbox{if (i) is satisfied}\\
|\nu|=1 \hbox{ and }\langle\xi_1,\nu\rangle=0&\hbox{if (ii) is satisfied}\\
|\nu|=1 \hbox{ and }\langle\xi_2,\nu\rangle=0&\hbox{if (iii) is satisfied}\end{array}\right.
$$
($\nu_1$, $\nu_2$, $\nu_3$ are the components of the vector $\nu$). Then, 
$$
\xi+\nabla\phi(x)=\left\{
\begin{array}{ll}
(\xi_1-\nu\mid \xi_2+\nu)&\hbox{if }x\in{\rm int}(\Delta_1)\\
(\xi_1-\nu\mid \xi_2-\nu)&\hbox{if }x\in{\rm int}(\Delta_2)\\
(\xi_1+\nu\mid \xi_2-\nu)&\hbox{if }x\in{\rm int}(\Delta_3)\\
(\xi_1+\nu\mid \xi_2+\nu)&\hbox{if }x\in{\rm int}(\Delta_4)
\end{array}
\right.
$$
(where ${\rm int}(E)$ denotes the interior of the set $E$). Taking Proposition \ref{propAA1}(i) into account, it follows that 
\begin{eqnarray}\label{Z_1}
\Z W_0(\xi)&\leq&{1\over 4}\Big(W_0(\xi_1-\nu\mid \xi_2+\nu)+W_0(\xi_1-\nu\mid \xi_2-\nu)\\
&&+\ W_0(\xi_1+\nu\mid \xi_2-\nu)+W_0(\xi_1+\nu\mid \xi_2+\nu)\Big).\nonumber
\end{eqnarray}
But
$
|(\xi_1-\nu)\land(\xi_2+\nu)|^2=|\xi_1\land \xi_2+(\xi_1+\xi_2)\land\nu|^2
=|\xi_1\land \xi_2|^2+|(\xi_1+\xi_2)\land\nu|^2
\geq|(\xi_1+\xi_2)\land\nu|^2,
$
and so 
$$
|(\xi_1+\nu)\land(\xi_2-\nu)|\geq |(\xi_1+\xi_2)\land\nu|=|\xi_1+\xi_2|.
$$
Similarly, we obtain: 
\begin{itemize}
\item[] $|(\xi_1-\nu)\land(\xi_2-\nu)|\geq |\xi_1-\xi_2|$;
\item[] $|(\xi_1+\nu)\land(\xi_2-\nu)|\geq |\xi_1+\xi_2|$;
\item[] $|(\xi_1+\nu)\land(\xi_2+\nu)|\geq |\xi_1-\xi_2|$.
\end{itemize}
Thus, $|(\xi_1-\nu)\land(\xi_2+\nu)|\geq\delta$, $|(\xi_1-\nu)\land(\xi_2-\nu)|\geq\delta$, $|(\xi_1+\nu)\land(\xi_2-\nu)|\geq\delta$ and $|(\xi_1+\nu)\land(\xi_2+\nu)|\geq\delta$, because $\min\{|\xi_1+\xi_2|,|\xi_1-\xi_2|\}\geq\delta$. Using ($\overline{\rm C}_2$) it follows that 
\begin{eqnarray*}
W_0(\xi_1-\nu\mid \xi_2+\nu)&\leq&c_\delta(1+|(\xi_1-\nu\mid \xi_2+\nu)|^p)\\
&\leq&c_\delta 2^p(1+|(\xi_1\mid \xi_2)|^p+|(-\nu\mid \nu)|^p)\\
&\leq&c_\delta 2^{2p+1}(1+|\xi|^p).
\end{eqnarray*}
In the same manner, we have: 
\begin{itemize}
\item[] $W_0(\xi_1-\nu\mid \xi_2-\nu)\leq c_\delta 2^{2p+1}(1+|\xi|^p)$;
\item[] $W_0(\xi_1+\nu\mid \xi_2-\nu)\leq c_\delta 2^{2p+1}(1+|\xi|^p)$;
\item[] $W_0(\xi_1+\nu\mid \xi_2+\nu)\leq c_\delta 2^{2p+1}(1+|\xi|^p)$,
\end{itemize}
and, from (\ref{Z_1}), we conclude that
$
\Z W_0(\xi)\leq c_\delta 2^{2p+1}(1+|\xi|^p).
$
\end{proof}

%%%%%%%%%%

\noindent{\em Proof of Proposition {\rm\ref{theorAA2}}. }Let $\xi=(\xi_1\mid \xi_2)\in\MM^{3\times 2}$. Then, one the four possibilities holds:
\begin{itemize}
\item[(i)] $|\xi_1\land\xi_2|\not=0$;
\item[(ii)] $|\xi_1\land\xi_2|=0$ with $\xi_1=\xi_2=0$;
\item[(iii)] $|\xi_1\land\xi_2|=0$ with $\xi_1\not=0$;
\item[(iv)] $|\xi_1\land\xi_2|=0$ with $\xi_2\not=0$.
\end{itemize}
For each $t\in\RR$, define $\varphi_t\in\Aff_0(Y;\RR)$ by 
$$
\varphi_t(x_1,x_2):=\left\{
\begin{array}{ll}
tx_2&\hbox{if }(x_1,x_2)\in\Delta_1\\
t(1-x_1)&\hbox{if }(x_1,x_2)\in\Delta_2\\
t(1-x_2)&\hbox{if }(x_1,x_2)\in\Delta_3\\
tx_1&\hbox{if }(x_1,x_2)\in\Delta_4
\end{array}
\right.
$$
with
\begin{itemize} 
\item[]$
\Delta_1:=\big\{(x_1,x_2)\in Y:x_2\leq x_1\leq -x_2+1\big\};
$
\item[]$
\Delta_2:=\big\{(x_1,x_2)\in Y:-x_1+1\leq x_2\leq x_1\big\};
$
\item[]$
\Delta_3:=\big\{(x_1,x_2)\in Y:-x_2+1\leq x_1\leq x_2\big\};
$
\item[]$
\Delta_4:=\big\{(x_1,x_2)\in Y:x_1\leq x_2\leq -x_1+1\big\}.
$
\end{itemize}
Consider $\phi\in\Aff_0(Y;\RR^3)$ given by 
$$
\phi:=\big(\varphi_{\nu_1},\varphi_{\nu_2},\varphi_{\nu_3}\big)\hbox{ with }
\left\{
\begin{array}{ll}
\nu={(\xi_1\land \xi_2)\over|\xi_1\land \xi_2|}&\hbox{if (i) is satisfied}\\  
|\nu|=1&\hbox{if (ii) is satisfied}\\
|\nu|=1\hbox{ and }\langle\xi_1,\nu\rangle=0& \hbox{if (iii) is satisfied}\\
|\nu|=1\hbox{ and }\langle\xi_2,\nu\rangle=0& \hbox{if (iv) is satisfied}
\end{array}
\right.
$$
($\nu_1$, $\nu_2$, $\nu_3$ are the components of the vector $\nu$). Then, 
$$
\xi+\nabla\phi(x)=\left\{
\begin{array}{ll}
(\xi_1\mid \xi_2+\nu)&\hbox{if }x\in{\rm int}(\Delta_1)\\
(\xi_1-\nu\mid \xi_2)&\hbox{if }x\in{\rm int}(\Delta_2)\\
(\xi_1\mid \xi_2-\nu)&\hbox{if }x\in{\rm int}(\Delta_3)\\
(\xi_1+\nu\mid \xi_2)&\hbox{if }x\in{\rm int}(\Delta_4)
\end{array}
\right.
$$
(where ${\rm int}(E)$ denotes the interior of the set $E$). Taking Proposition \ref{propAA1}(ii) into account, it follows that
\begin{eqnarray}\label{ZzZ}
\Z W_0(\xi)&\leq&{1\over 4}\Big(\Z W_0(\xi_1\mid \xi_2+\nu)+\Z W_0(\xi_1-\nu\mid \xi_2)\\
&&+\ \Z W_0(\xi_1\mid \xi_2-\nu)+\Z W_0(\xi_1+\nu\mid \xi_2)\Big).\nonumber
\end{eqnarray}
But
$
|\xi_1+(\xi_2+\nu)|^2=|(\xi_1+\xi_2)+\nu|^2
=|\xi_1+\xi_2|^2+|\nu|^2
=|\xi_1+\xi_2|^2+1
\geq 1,
$
hence
$
|\xi_1+(\xi_2+\nu)|\geq 1.
$
Similarly, we obtain
$
|\xi_1-(\xi_2+\nu)|\geq 1,
$
and so
$$
\min\{|\xi_1+(\xi_2+\nu)|,|\xi_1-(\xi_2+\nu)|\}\geq 1.
$$
In the same manner, we have: 
\begin{itemize}
\item[] $\min\{|(\xi_1-\nu)+\xi_2|,|(\xi_1-\nu)-\xi_2|\}\geq 1$;
\item[] $\min\{|\xi_1+(\xi_2-\nu)|,|\xi_1-(\xi_2-\nu)|\}\geq 1$;
\item[] $\min\{|(\xi_1+\nu)+\xi_2|,|(\xi_1+\nu)-\xi_2|\}\geq 1$. 
\end{itemize}
Using Lemma \ref{lemmaA3} it follows that
\begin{eqnarray*}
\Z W_0(\xi_1\mid \xi_2+\nu)&\leq& r_1(1+|(\xi_1\mid \xi_2+\nu)|^p)\\
&\leq&r_1 2^p(1+|(\xi_1\mid\xi_2)|^p+|(0\mid\nu)|^p)\\
&\leq&r_12^{p+1}(1+|\xi|^p).
\end{eqnarray*}
Similarly, we obtain:
\begin{itemize}
\item[] $\Z W_0(\xi_1-\nu\mid \xi_2)\leq r_12^{p+1}(1+|\xi|^p)$;
\item[] $\Z W_0(\xi_1\mid \xi_2-\nu)\leq r_12^{p+1}(1+|\xi|^p)$;
\item[] $\Z W_0(\xi_1+\nu\mid \xi_2)\leq r_12^{p+1}(1+|\xi|^p)$,
\end{itemize} 
and, from (\ref{ZzZ}), we conclude  that 
$
\Z W_0(\xi)\leq r_12^{p+1}(1+|\xi|^p).
$
\hfill$\square$
 
 %%%%%%%%%%%
 
\medskip

The next proposition will be used in the  proof of Theorem \ref{anzman}.

\begin{proposition}\label{corA4}
If {\rm($\overline{{\rm C}}_2$)} holds then $\Z W_0=\mathcal{Q}W_0=\mathcal{Q}[\Z W_0]$. 
\end{proposition}
\begin{proof}
By Proposition \ref{theorAA2}, $\Z W_0(\xi)\leq c(1+|\xi|^p)$ for all $\xi\in\MM^{3\times 2}$ and some $c>0$. Then $\Z W_0$ is finite, and so $\Z W_0$ is continuous by Proposition \ref{propAA1}(iii). Recall the  (classical) theorem:
\begin{theorem}[Dacorogna \cite{dacorogna}]\label{DQF} If $f:\MM^{3\times 2}\to[0,+\infty]$ is finite and continuous then $\Z f=\mathcal{Q}f$. 
\end{theorem} 
\noindent By Theorem \ref{DQF} we have $\Z[\Z W_0]=\mathcal{Q}[\Z W_0]$. But $\Z[\Z W_0]=\Z W_0$ by Proposition \ref{propAA1}(ii), hence $\Z W_0=\mathcal{Q}[\Z W_0]$. Thus $\Z W_0$ is quasiconvex and $\Z W_0\leq W_0$. On the other hand, noticing that $\Z g=g$ whenever $g$ is quasiconvex, we see that if $g$ is quasiconvex  and $g\leq W_0$ then $g\leq \Z W_0$. According to Definition \ref{DEFofQUasiRankoNEConvexityandEnvelOPe}(ii), it follows that $\Z W_0=\mathcal{Q}W_0$.
 \end{proof}

                      %%%%%% Proof of Theorem \ref{anzman} %%%%%%%%%

\subsection{Proof of Theorems \ref{add} and \ref{anzman}} We begin by proving Proposition \ref{fund_prop}  which will play an essential role in the proof of Theorems \ref{add} and \ref{anzman}. 
\begin{proposition}\label{fund_prop}
$
\overline{\J}=\mathcal{J}
$
with $\mathcal{J}:W^{1,p}(\Sigma;\RR^3)\to[0,+\infty]$ given by
$$
\mathcal{J}(v):=\inf\left\{\liminf_{n\to+\infty}\int_\Sigma \Z W_0(\nabla v_n(x))dx:\Aff(\Sigma;\RR^3)\ni v_n\to v\hbox{ in }L^{p}(\Sigma;\RR^3)\right\}.
$$
\end{proposition}
To prove Proposition \ref{fund_prop} we need the following lemma.

\begin{lemma}\label{fund_lem}
If $v\in\Aff(\Sigma;\RR^3)$ then
\begin{equation}\label{step1}
\overline{\J}(v)\leq \int_\Sigma \Z W_0(\nabla v(x))dx.
\end{equation}
\end{lemma}
\begin{proof}
Let $v\in\Aff(\Sigma;\RR^3)$. By definition, there exists a finite family $(V_i)_{i\in I}$ of open disjoint subsets of $\Sigma$ such that $|\partial V_i|=0$ for all $i\in I$, $|\Sigma\setminus\cup_{i\in I}V_i|=0$ and, for every $i\in I$, $\nabla v(x)=\xi_i$ in $V_i$ with $\xi_i\in\MM^{3\times 2}$. Given any $\delta>0$ and any $i\in I$, we consider $\phi_i\in \Aff_0(Y;\RR^3)$ such that
\begin{equation}\label{Zk}
\int_Y W_0(\xi_i+\nabla\phi_i(y))dy\leq\Z W_0(\xi_i)+{\delta\over|\Sigma|}.
\end{equation}
Fix any integer $n\geq 1$. By Vitali's covering theorem, there exists a finite or countable family $(a_{i,j}+\alpha_{i,j}Y)_{j\in J_{i}}$ of disjoint subsets of $V_i$, where $a_{i,j}\in\RR^2$ and $0<\alpha_{i,j}<{1\over n}$, such that
$
\big|V_i\setminus\cup_{j\in J_{i}}(a_{i,j}+\alpha_{i,j}Y)\big|=0
$
 (and so $\sum_{j\in J_i}\alpha_{i,j}^2=|V_i|$). Define $\psi_n:\Sigma\to\RR^3$ by
$$
\psi_n(x):=
\alpha_{i,j}\phi_{i}\left({x-a_{i,j}\over \alpha_{i,j}}\right)\hbox{ if }x\in a_{i,j}+\alpha_{i,j}Y.
$$
Since $\phi_i\in\Aff_0(Y;\RR^3)$, there exists a finite family $(Y_{i,l})_{l\in L_i}$ of open disjoint subsets of $Y$ such that $|\partial Y_{i,l}|=0$ for all $l\in L_i$, $|Y\setminus\cup_{l\in L_i}Y_{i,l}|=0$ and, for every $l\in L_i$, $\nabla\phi_i(y)=\zeta_{i,l}$ in $Y_{i,l}$ with $\zeta_{i,l}\in\MM^{3\times 2}$. Set 
$
U_{i,l,n}:=\cup_{j\in J_i}a_{i,j}+\alpha_{i,j}Y_{i,l},
$
then $|\partial U_{i,l,n}|=0$ for all $i\in I$ and all $l\in L_i$, $|\Sigma\setminus\cup_{i\in I}\cup_{l\in L_i}U_{i,l,n}|=0$ and, for every $i\in I$ and every $l\in L_i$, $\nabla\psi_n(x)=\zeta_{i,l}$ in $U_{i,l,n}$, and so $\psi_n\in\Aff_0(\Sigma;\RR^3)$. On the other hand, $\|\psi_n\|_{L^\infty(\Sigma;\RR^3)}\leq {1\over n}\max_{i\in I}\|\phi_i\|_{L^\infty(Y;\RR^3)}$ and $\|\nabla\psi_n\|_{L^\infty(\Sigma;\MM^{3\times 2})}\leq\max_{i\in I}\|\nabla\phi_i\|_{L^\infty(Y;\MM^{3\times 2})}$, hence (up to a subsequence) $\psi_n\stackrel{*}{\wto}0$ in $W^{1,\infty}(\Sigma;\RR^3)$, where ``$\stackrel{*}{\wto}$" denotes the weak$^*$ convergence in $W^{1,\infty}(\Sigma;\RR^3)$. Consequently, $\psi_n\wto 0$ in $W^{1,p}(\Sigma;\RR^3)$, and so (up to a subsequence) $\psi_n\to 0$ in $L^p(\Sigma;\RR^3)$. Moreover, 
\begin{eqnarray*}
\int_\Sigma W_0\left(\nabla v(x)+\nabla\psi_n(x)\right)dx&=&\sum_{i\in I}\int_{V_i} W_0\left(\xi_i+\nabla\psi_n(x)\right)dx\\
&=&\sum_{i\in I}\sum_{j\in J_i}\alpha_{i,j}^2\int_{Y}W_0\left(\xi_i+\nabla\phi_{i}(y)\right)dy\\
&=&\sum_{i\in I}|V_i|\int_{Y}W_0\left(\xi_i+\nabla\phi_{i}(y)\right)dy.
\end{eqnarray*}
As $v+\psi_n\in\Aff(\Sigma;\RR^3)$ and $v+\psi_n\to v$ in $L^{p}(\Sigma;\RR^3)$, from (\ref{Zk}) we deduce that
\begin{eqnarray*}
\overline{\J}(v)\leq\liminf_{n\to+\infty}\int_\Sigma W_0\left(\nabla v(x)+\nabla\psi_n(x)\right)dx&\leq&\sum_{i\in I}|V_i|\Z W_0(\xi_i)+\delta\\
&=&\int_\Sigma\Z W_0(\nabla v(x))dx + \delta,
\end{eqnarray*}
and (\ref{step1}) follows.
\end{proof}
\begin{remark}\label{CPAFremark2}
As the proof of Lemma \ref{fund_lem} requires Vitali's covering theorem, if we consider  ``$\Aff^{ET}$" (with ``$\Aff^{ET}$" defined in Remark \ref{EkelandTemamDensity}) instead of ``$\Aff$", Lemma \ref{fund_lem} seems to be false. However, Lemma \ref{fund_lem} remains valid if we replace ``$\Aff$" by ``$\Aff^V$" (with ``$\Aff^V$" defined in Remark \ref{CPAFremark}). 
\end{remark}
\noindent{\em Proof of Proposition {\rm \ref{fund_prop}}.  }Clearly $\mathcal{J}\leq\overline{\J}$. We are thus reduced to prove that 
\begin{equation}\label{inequality}
\overline{\J}\leq \mathcal{J}.
\end{equation}
Fix any $v\in W^{1,p}(\Sigma;\RR^3)$ and any sequence $v_n\to v$ in $L^p(\Sigma;\RR^3)$ with $v_n\in\Aff(\Sigma;\RR^3)$. Using Lemma \ref{fund_lem} we have
$
\overline{\J}(v_n)\leq\int_\Sigma \Z W_0(\nabla v_n(x))dx
$
for all $n\geq 1$. Thus,
$$
\overline{\J}(v)\leq\liminf_{n\to+\infty}\overline{\J}(v_n)\leq\liminf_{n\to+\infty}\int_\Sigma \Z W_0(\nabla v_n(x))dx,
$$
and (\ref{inequality}) follows.\hfill $\square$

\medskip

\noindent{\em Proof of Theorem {\rm\ref{add}}. }By Proposition \ref{theorAA2}, $\Z W_0(\xi)\leq c(1+|\xi|^p)$ for all $\xi\in\MM^{3\times 2}$ and some $c>0$. Then $\Z W_0$ is finite, and so $\Z W_0$ is continuous by Proposition \ref{propAA1}(iii). As $\Aff(\Sigma;\RR^3)$ is strongly dense in $W^{1,p}(\Sigma;\RR^3)$, we deduce that for every $v\in W^{1,p}(\Sigma;\RR^3)$,
$$
\mathcal{J}(v)=\inf\left\{\liminf_{n\to+\infty}\int_\Sigma \Z W_0(\nabla v_n(x))dx:W^{1,p}(\Sigma;\RR^3)\ni v_n\to v\hbox{ in }L^{p}(\Sigma;\RR^3)\right\},
$$
and so $\mathcal{J}\leq \mathcal{I}$. But $\mathcal{I}\leq\overline{\J}$ and $\overline{\J}=\mathcal{J}$ by Proposition \ref{fund_prop}, hence $\overline{\J}=\mathcal{I}$.\hfill$\square$

\medskip

\noindent{\em Proof of Theorem {\rm\ref{anzman}}. }An analysis similar to that of the proof of Theorem \ref{add} shows that $\Z W_0$ is continuous, $\Z W_0(\xi)\leq c(1+|\xi|^p)$ for all $\xi\in\MM^{3\times 2}$ and some $c>0$, and 
$$
\overline{\J}(v)=\inf\left\{\liminf_{n\to+\infty}\int_\Sigma \Z W_0(\nabla v_n(x))dx:W^{1,p}(\Sigma;\RR^3)\ni v_n\to v\hbox{ in }L^{p}(\Sigma;\RR^3)\right\}.
$$
Recall the  (classical) integral representation theorem:
\begin{theorem}[Dacorogna \cite{dacorogna}]\label{DIRT}
Let $f:\MM^{3\times 2}\to[0,+\infty]$ be a Borel measurable function and let $\mathcal{F}:W^{1,p}(\Sigma;\RR^3)\to[0,+\infty]$ be defined by
$$
\mathcal{F}(v):=\inf\left\{\liminf_{n\to+\infty}\int_\Sigma f(\nabla v_n(x))dx:W^{1,p}(\Sigma;\RR^3)\ni v_n\to v\hbox{ in }L^{p}(\Sigma;\RR^3)\right\}.
$$ 
If $f$ is continuous and $C|\xi|^p\leq f(\xi)\leq c(1+|\xi|^p)$ for all $\xi\in\MM^{3\times 2}$ and some $c,C>0$, then for every $v\in W^{1,p}(\Sigma;\RR^3)$,
$$
\mathcal{F}(v)=\int_\Sigma \mathcal{Q}f(\nabla v(x))dx.
$$
\end{theorem}
\noindent Noticing that $\Z W_0$ is coercive, from Theorem \ref{DIRT} it follows that for every $v\in W^{1,p}(\Sigma;\RR^3)$,
$$
\overline{\J}(v)=\int_\Sigma\mathcal{Q}[\Z W_0](\nabla v(x))dx.
$$
Moreover, $\mathcal{Q}[\Z W_0]=\mathcal{Q} W_0$ by Proposition \ref{corA4}, and the proof is complete.\hfill$\square$

\bigskip

\noindent{\em Acknowledgments. }The authors wish to thank an anonymous referee for interesting suggestions and for giving a simpler proof of Lemma \ref{Lemma3}.

%%%%%%%%%%%%%%%% REFERENCES %%%%%%%%%%%%%%%%%%%

\end{document}